\documentclass[12pt, a4paper, oneside]{article}
\usepackage[top=1in, bottom=1.1in, left=0.9in, right=0.5in]{geometry}
\usepackage[parfill]{parskip}
\usepackage{amsmath}
\usepackage{amsthm}
\usepackage{graphicx}
\usepackage{amssymb}
\usepackage{epstopdf}
\usepackage{setspace}
\usepackage{authblk}
\usepackage{enumerate}
\usepackage{lmodern}
\usepackage[hypcap]{caption}
\usepackage[english]{babel}
\usepackage{fancyhdr}
\addto\captionsenglish{}
\addto\captionsenglish{}
\addto\captionsenglish{}
\addto\captionsenglish{}
\addto\captionsenglish{}

\newlength\longest

\newcommand{\R}{{\mathbb R}}

\begin{document}

\newtheorem{theorem}{Teorema}[section]
\newtheorem{lemma}[theorem]{Lemma}
\newtheorem{corollary}[theorem]{Corollary}
\newtheorem{proposition}[theorem]{Proposition}
\newtheorem{conjecture}[theorem]{Conjecture}
\newtheorem{problem}[theorem]{Problem}
\newtheorem{claim}[theorem]{Claim}
\theoremstyle{definition}
\newtheorem{assumption}[theorem]{Assumption}
\newtheorem{remark}[theorem]{Remark}
\newtheorem{definition}[theorem]{Definition}
\newtheorem{example}[theorem]{Example}
\theoremstyle{remark}
\newtheorem{notation}{Notasi}
\renewcommand{\thenotation}{}

\title{Inclusion Properties of Weighted Weak Orlicz Spaces}
\author{ Al Azhary Masta${}^1$, Ifronika${}^{2}$, Muhammad Taqiyuddin${}^3$}

\affil{${}^{1,2}$Department of Mathematics, Institut Teknologi Bandung\\ Jl. Ganesha no. 10, Bandung}

\affil{${}^{1,3}$\emph{Permanent Address}: Department of Mathematics Education, Universitas Pendidikan Indonesia, Jl. Dr. Setiabudi 229, Bandung 40154\\

E-mail: ${}^{1}$alazhari.masta@upi.edu, ${}^{2}$ifronika@gmail.com,
${}^{3}$taqi94@hotmail.com}
\date{}

\maketitle

\begin{abstract}
In this paper we discuss the structure of weighted weak Lebesgue spaces and weighted weak Orlicz spaces on $\mathbb{R}^n$. First, we present sufficient and necessary conditions for inclusion relation between weighted weak Lebesgue spaces. Next, we also obtain similar results on weighted weak Orlicz spaces. One of the keys to prove our results is to use the norm of the characteristic functions of the balls in $\mathbb{R}^n$.
\bigskip

\end{abstract}

\section{Introduction}

Orlicz spaces as generalization of Lebesgue spaces were introduced by Z. W. Birnbaum and W. Orlicz in 1931 (see \cite{Kufner,Luxemburg,Lech,Rao}). Many authors gave more attention in the study of Orlicz spaces (see \cite{Yong,Krasnosel'skii,Kufner,Christian,Ning,Luxemburg,Lech,Orlicz,Rao,Welland,Xueying}). In particular, Maligranda \cite{Lech} discussed about inclusion properties of Orlicz spaces. In 2016, Masta\textit{ et al.} \cite{Masta1} obtained inclusion relations between Orlicz spaces and between weak Orlicz spaces using different technique.

On the other hand, Osan\c{c}liol \cite{Alen} have proved sufficient and necessary conditions for inclusion relation between two weighted Lebesgue spaces. Furthermore, he also obtained the sufficient and necessary conditions for inclusion relation between two weighted (strong) Orlicz spaces.

Motivated by these results, we are interested in studying the inclusion properties of weighted weak Lebesgue spaces and weighted weak Orlicz spaces. In particular, we shall give sufficient and necessary conditions for inclusion relations between weighted Lebesgue spaces and between weighted weak Orlicz spaces.

For that purpose, we introduce some definitions. Let $U$ be the set of all functions $u : \mathbb{R}^n \rightarrow ( 0, \infty )$ such that $u(x+y)\leq u(x)\cdot u(y)$ for all $ x, y \in \mathbb{R}^n$. Let $u_1, u_2$ be the elements of $U,$ we denote $u_1 \preceq u_2$ if there exists a constant $C > 0$ such that $u_1(x)  \leq C u_2(x)$ for all $ x \in \mathbb{R}^n$. 

Now, let us give the definition of weighted weak Lebesgue spaces and weighted weak Orlicz spaces.  For $ 1 \leq p < \infty$ and $ u: \R^n \rightarrow (0,\infty)$, the weighted weak Lebesgue space $wL_{p}^u(\mathbb{R}^n)$ is the set of all measurable function $f : \mathbb{R}^n \rightarrow \mathbb{R}$  such that $$\| f \|_{wL_{p}^u(\mathbb{R}^n)}:= \mathop {\sup }\limits_{t > 0} t | \{ x \in\mathbb{R}^n : |u(x)f(x)| > t \} |^{\frac{1}{p}} < \infty .$$ 

Let $\Phi: [ 0, \infty ) \rightarrow [ 0, \infty ) $ be a Young function [that is,  $\Phi$ is convex, $\lim \limits_{t\to 0}\Phi(t)=0=\Phi(0)$, left-continuous and $ \lim \limits_{t\to\infty} \Phi(t) = \infty$] and $u : \mathbb{R}^n \rightarrow ( 0, \infty )$, the weighted weak Orlicz space $wL_{\Phi}^u(\mathbb{R}^n)$ is the set of all measurable function $f : \mathbb{R}^n \rightarrow \mathbb{R}$  such that
$$
\| f \|_{wL_{\Phi}^u(\mathbb{R}^n)} := \inf \left\{  {b>0:
	\mathop {\sup }\limits_{t > 0} \Phi(t) \Bigl| \Bigl\{ x \in \mathbb{R}^n : \frac{|u(x)f(x)|}{b} > t \Bigr\}
	\Bigr| \leq1}\right\} < \infty.
$$

Note that if there exists $C>0$ such that $u_1 (x) \leq Cu_2 (x)$ for every $x \in \mathbb{R}^n$, then $$\| f \|_{wL_{\Phi}^{u_1}(\mathbb{R}^n)} \leq C \| f \|_{wL_{\Phi}^{u_2}(\mathbb{R}^n)}.$$  

The space $wL_{\Phi}^u(\mathbb{R}^n)$ is quasi-Banach spaces equipped with the quasi-norm $\| f \|_{wL_{\Phi}^u(\mathbb{R}^n)}$. If $u(x)=1$ for every $ x \in \R^n$, then $wL_{\Phi}^u(\mathbb{R}^n)$ is weak Orlicz space $wL_{\Phi}(\mathbb{R}^n)$. Meanwhile, for $\Phi(t)=t^p$ ($ 1 \leq p < \infty$), we have $wL_{\Phi}^u(\mathbb{R}^n) = wL_{p}^u(\mathbb{R}^n)$. 
  
To achieve our purpose, we will use the similar methods in \cite{Gunawan, Masta1, Masta2, Alen} which pay attention to the characteristic functions of open balls in $\R^n$. Next, we recall some lemmas which will be used later in next section.

\medskip

\begin{lemma}\label{lemma:1.1} \cite{Nakai1}
Suppose that $\Phi$ is a Young function and $ \Phi^{-1}(s):=\inf \{r \geq 0 :
\Phi (r) > s \}$. We have

{\parindent=0cm
{\rm (1)} $\Phi^{-1}(0) = 0$.

{\rm (2)} $ \Phi^{-1}(s_1) \leq \Phi^{-1}(s_2)$ for  $s_1 \leq s_2$.

{\rm (3)} $\Phi (\Phi^{-1}(s)) \leq s \leq \Phi^{-1}(\Phi(s))$ for $0 \leq s <
\infty$.

\par}

\end{lemma}

\bigskip

\begin{lemma} \cite{Masta3}\label{lemma:1.2}
Let $\Phi_1, \Phi_2$ be Young functions. For any $s>0$, if there exist constants $C_1, C_2 > 0$ such that $\Phi^{-1}_2(s) \leq C_1 \Phi^{-1}_1(C_2s)$, then we have  $\Phi_1 (\frac{t}{C_1})  \leq C_2 \Phi_2(t)$ for $ t = \Phi^{-1}_2(s).$
\end{lemma}

\bigskip

\bigskip
	\begin{lemma}\label{lemma:1.3}
 Let $u \in U$ and $ 1 \leq p < \infty$. If $L_{x}f$ is translation function, i.e $L_{x}f(y):= f(y-x)$ for every $ x \in \mathbb{R}^n$, then:
	{\parindent=0cm
		
	{\rm (1)}  For all $ f \in wL_{p}^{u}(\mathbb{R}^n)$ and for all $ x \in \mathbb{R}^n$, we have $L_{x}f \in wL_{p}^{u}(\mathbb{R}^n)$ and $$\|L_{x}f\|_{wL_{p}^{u}(\mathbb{R}^n)} \leq u(x)\|f\|_{wL_{p}^{u}(\mathbb{R}^n)}.$$  
		
	{\rm (2)} If $ f \in wL_{p}^{u}(\mathbb{R}^n)$ and $f\neq0$, then there exist a constant $C > 0$ (depends on $f$) such that $$\frac{u(x)}{C} \leq \|L_{x}f\|_{wL_{p}^{u}(\mathbb{R}^n)} \leq C u(x).$$
		
		\par}
\end{lemma}

In this paper, the letter $C$ will be used for constants that may change from line to line, while constants with subscripts, such as $C_{1}, C_{2}$, do not change in different lines.

\section{Results}

First, we present sufficient and necessary conditions for inclusion properties of weighted weak Lebesgue spaces in the following theorem.

\medskip

\begin{theorem}\label{theorem:2.1}
Let $ 1 \leq p < \infty$ and $ u_1, u_2: \mathbb{R}^n \rightarrow (0,\infty)$. Then the following statements are equivalent:
	
{\parindent=0cm

{\rm (1)}  $u_1 \preceq u_2$.
		
{\rm (2)} $wL_{p}^{u_2}(\mathbb{R}^n) \subseteq wL_{p}^{u_1}(\mathbb{R}^n)$.

{\rm (3)} There exist a constant $C>0$ such that $$\|f\|_{wL_{p}^{u_1}(\mathbb{R}^n)} \leq C \|f\|_{wL_{p}^{u_2}(\mathbb{R}^n)},$$ for every $ f \in wL_{p}^{u_2}(\mathbb{R}^n)$.
		
\par}
\end{theorem}

Proof.

Assume that (1) holds, then there exist a constant $C>0$ such that $ u_1(x) \leq C u_2(x)$ for every $ x \in \R^n$. Let  $ f$ be an element of $ wL_{p}^{u_2}(\mathbb{R}^n)$.   Now, take an arbitrary $ t > 0$, observe that (by setting $t_1 = \frac{t}{C}$)

\begin{align*}
t |\{x \in \R^n : |u_1(x)f(x)| > t \}|^{\frac{1}{p}}  & \leq   t |\{x \in \R^n : |Cu_2(x)f(x)| > t \}|^{\frac{1}{p}} \\
&=   t \Bigl| \Bigl\{x \in \R^n : |u_2(x)f(x)| > \frac{t}{C} \Bigr\}\Bigr|^{\frac{1}{p}} \\
& = Ct_1 |\{x \in \R^n : |u_2(x)f(x)| > t_1 \}|^{\frac{1}{p}}\\
& \leq C \|  f \|_{wL_{p}^{u_2}(\R^n)}.
\end{align*}

Since $ t > 0$ is arbitrary, we have $\|  f \|_{wL_{p}^{u_1}(\R^n)} \leq C \|  f \|_{wL_{p}^{u_2}(\R^n)}$.  

 Next, since $(wL_{p}^{u_1}(\R^n), wL_{p}^{u_2}(\R^n))$ is a quasi-Banach pair, as mentioned in [10, Appendix G], we are aware that [6, Lemma 3.3] still holds for quasi-Banach spaces, so we have (2) and (3) are equivalent. Then we remain to show that (3) implies (1).
Assume that (3) holds. By Lemma \ref{lemma:1.3}, we have

$$\frac{u_{1}(x)}{C}\leq \|L_{x}f\|_{wL_{p}^{u_1}(\mathbb{R}^n)} \leq C \|L_{x}f\|_{wL_{p}^{u_2}(\mathbb{R}^n)}\leq u_2(x).$$

So, we can conclude that $u_1 \preceq u_2$. \qed

Next, we will investigate the inclusion properties of weighted weak Orlicz spaces. For getting a result, we give attention to estimate the norm of the characteristic function of open ball in the following lemma.

\bigskip

\begin{lemma}\label{lemma:2.2}\cite{Ning,Masta1}
	Let $\Phi$ be a Young function, $a\in\R^n$, and $r >0$ be arbitrary. Then we have
	$ \left\|\frac{\chi_{B(a,r)}}{u} \right\|_{wL_{\Phi}^u(\mathbb{R}^n)} = \|\chi_{B(a,r)}\|_{wL_\Phi(\mathbb{R}^n)} = \frac{1}{\Phi^{-1}\bigl( \frac{1}{|B(a,r)|} \bigr)}$ where $|B(a,r)|$ denotes the volume of open ball $B(a,r)$ centered at $a\in\R^n$ with radius $r >0$.
\end{lemma}

\medskip

Now we come to the inclusion property of weighted weak Orlicz spaces $wL_{\Phi_1}^u(\mathbb{R}^n)$ and $wL_{\Phi_2}^u(\mathbb{R}^n)$ with respect to Young functions $\Phi_1, \Phi_2$. Given two Young functions $\Phi_1, \Phi_2$, we write $\Phi_1 \prec \Phi_2$ if there exists a constant $C > 0$ such that $\Phi_1(t) \leq \Phi_2(Ct)$ for all $ t > 0$.    

\medskip
  
\begin{theorem}\label{theorem:2.3}
	Let $\Phi_1, \Phi_2$ be Young functions and $ u: \mathbb{R}^n \rightarrow (0,\infty)$. Then the following statements are equivalent:
	
	{\parindent=0cm
		{\rm (1)} $\Phi_1 \prec \Phi_2$.
		
		{\rm (2)} $wL_{\Phi_2}^u(\mathbb{R}^n) \subseteq wL_{\Phi_1}^u(\mathbb{R}^n)$.
		
		{\rm (3)} There exist a constant $C > 0$ such that $\| f \|_{wL_{\Phi_1}^u(\mathbb{R}^n)} \leq C\| f \|_{wL_{\Phi_2}^u(\mathbb{R}^n)}$, for every $ f \in wL_{\Phi_2}^u(\mathbb{R}^n)$.
		
		\par}
\end{theorem}

\noindent{\it Proof}.
Assume that (1) holds. Let $ f \in wL_{\Phi_2}^u(\mathbb{R}^n)$,
$$A_{\Phi_1,u} = \Bigl\{ b>0: \mathop {\sup }\limits_{t > 0} \Phi_1(t) \left|
	\Bigl\{ x \in \R^n : \frac{|u(x)f(x)|}{b} > t \Bigr\}\right|  \leq1 \Bigr\}
$$
and
\begin{align*}
A_{\Phi_2,u} & = \Bigl\{  {b>0: \mathop {\sup }\limits_{t > 0} \Phi_2(Ct) \left|
	\Bigl\{ x \in \R^n : \frac{|u(x)f(x)|}{b} > t \Bigr\}\right| \leq1}\Bigr\} \\
& = \Bigl\{  {b>0: \mathop {\sup }\limits_{s > 0} \Phi_2(s) \Bigl|
	\Bigl\{ x \in \mathbb{R}^n : \frac{C|u(x)f(x)|}{b} > s \Bigr\} \Bigr| \leq1}\Bigr\}.
\end{align*}

Observe that, for arbitrary $ b \in A_{\Phi_2,u} $ and $ t > 0$, we have

\begin{align*}
\Phi_1(t) \Bigl| \Bigl\{x \in \R^n : \frac{|u(x)f(x)|}{b} > t \Bigr\} \Bigr|  & \leq   \Phi_2(Ct) \Bigl| \Bigr\{x \in \R^n : \frac{|u(x)f(x)|}{b} > t \Bigr\} \Bigr| \\
&\leq  \mathop {\sup }\limits_{s > 0} \Phi_2(s) \Bigl| \Bigl\{x \in \R^n : \frac{C|u(x)f(x)|}{b} > s \Bigr\} \Bigr| \\
&\leq 1.
\end{align*}

Since $t>0$ is arbitrary, we obtain
$$\| f \|_{wL_{\Phi_1}^u(\mathbb{R}^n)} := \inf A_{\Phi_1,u} \leq \inf A_{\Phi_2,u} = \| Cf \|_{wL_{\Phi_2}^u(\mathbb{R}^n)} := C \| f \|_{wL_{\Phi_2}^u(\mathbb{R}^n)},$$
which also proves that $wL_{\Phi_2}^u(\mathbb{R}^n) \subseteq wL_{\Phi_1}^u(\mathbb{R}^n)$.
Using a similar argument in the proof of Theorem \ref{theorem:2.1}, we have (2) and (3) are equivalent.

Assume now that (3) holds. By Lemma \ref{lemma:2.2}, we have
$$
\frac{1}{\Phi_1^{-1}\Bigl(\frac{1}{|B(a,r)|}\Bigr)} = \left\|  \frac{\chi_{B(a,r)}}{u} \right\|_{wL_{\Phi_1}^u(\mathbb{R}^n)}
\leq C \left\|\frac{\chi_{B(a,r)}}{u} \right\|_{wL_{\Phi_2}^u(\mathbb{R}^n)} =  \frac{C}{\Phi_2^{-1}\Bigl(\frac{1}{|B(a,r)|}\Bigr)}
$$
or $ C \Phi_1^{-1}(\frac{1}{|B(a,r)|}) \geq \Phi_2^{-1}(\frac{1}{|B(a,r)|})$,
for arbitrary $a\in\mathbb{R}^n$ and $r > 0$. By Lemma \ref{lemma:1.2}, we have
$$
\Phi_1\Bigl(\frac{C}{t_0}\Bigr) \leq  \Phi_2(t_0),
$$
for $t_0 = \Phi_2^{-1}( \frac{1}{|B(a,r)|})$. Since $a\in\mathbb{R}^n$ and $r > 0$ are arbitrary, we conclude that $\Phi_1(t) \leq \Phi_2(Ct)$ for every $t > 0$.\qed 

\medskip

\remark For $u(x)=1$, Theorem \ref{theorem:2.3} reduces to Theorem 3.3 in \cite{Masta1}. 

\medskip
Next, we also give the sufficient and necessary conditions for inclusion relation between weighted weak Orlicz spaces $wL_{\Phi_1}^{u_1}(\mathbb{R}^n)$ and $wL_{\Phi_2}^{u_2}(\mathbb{R}^n)$ with respect to Young functions $\Phi_1, \Phi_2$ and weights $u_1, u_2$. To get the result, we need some lemmas in the following.

\medskip

\begin{lemma}\label{lemma:2.3}
	Let $\Phi$ be a Young function. If $ f \in wL_{\Phi}^u(\mathbb{R}^n)$, then for arbitrary $\epsilon > 0$ we have $$ \mathop {\sup}\limits_{t > 0} \Phi(t) \left| \left\{ x \in \R^n : \frac{|u(x)f(x)|}{ \| f \|_{wL_{\Phi}^u(\mathbb{R}^n)} + \epsilon}  > t \right\} \right| \leq 1.$$
	
	Proof.
\end{lemma}
Let $f$ be an element of $wL_{\Phi}^u (\R^n)$. Take an arbitrary $\epsilon > 0$, then there exists $b_{\epsilon} > 0$ such that $b_{\epsilon} \leq \| f \|_{wL_{\Phi}^u(\mathbb{R}^n)} + \epsilon$ and $$ \mathop {\sup}\limits_{t > 0} \Phi(t) \left| \left\{ x \in \R^n : \frac{|u(x)f(x)|}{ b_{\epsilon}}  > t \right\} \right| \leq 1.$$ 

Because $\frac{|u(x)f(x)|}{b_{\epsilon}} \geq \frac{|u(x)f(x)|}{\| f \|_{wL_{\Phi}^u(\mathbb{R}^n)} + \epsilon}$, we have $$\Phi(t)\left| \left\{ x \in \R^n : \frac{|u(x)f(x)|}{ \| f \|_{wL_{\Phi}^u(\mathbb{R}^n)} + \epsilon}  > t \right\} \right|  \leq \Phi(t) \left| \left\{ x \in \R^n : \frac{|u(x)f(x)|}{ b_{\epsilon}}  > t \right\} \right| \leq 1$$ for every  $t > 0$.
Since $t > 0$ is arbitrary, we can conclude that$$ \mathop {\sup}\limits_{t > 0} \Phi(t) \left| \left\{ x \in \R^n : \frac{|u(x)f(x)|}{ \| f \|_{wL_{\Phi}^u(\mathbb{R}^n)} + \epsilon}  > t \right\} \right| \leq 1$$ for every $\epsilon > 0$. \qed

\medskip

\begin{lemma}\label{lemma:2.6}
Let $u$ be element of $U$. If $\Phi$ is a continuous Young function satisfying the $\bigtriangleup_{2}$ condition [that is, there exists $ K > 0$ such that $\Phi(2t) \leq K \Phi(t)$ for all $ t \geq 0$], then:
{\parindent=0cm

{\rm (1)}  For all $ f \in wL_{\Phi}^{u}(\mathbb{R}^n)$ and for all $ x \in \mathbb{R}^n$, we have $L_{x}f \in wL_{\Phi}^{u}(\mathbb{R}^n)$ and $$\|L_{x}f\|_{wL_{\Phi}^{u}(\mathbb{R}^n)} \leq u(x)\|f\|_{wL_{\Phi}^{u}(\mathbb{R}^n)}.$$

{\rm (2)} If $ f \in wL_{\Phi}^{u}(\mathbb{R}^n)$ and $f \neq 0$, then there exist a constant $ C > 0$ (depends on $f$) such that $$\frac{u(x)}{C} \leq \|L_{x}f\|_{wL_{\Phi}^{u}(\mathbb{R}^n)} \leq C u(x).$$

\par}
\end{lemma}
\noindent{\it Proof}.

{\parindent=0cm

{\rm (1)} Let $ f \in wL_{\Phi}^{u}(\mathbb{R}^n)$ and $L_{x}f(y)=f(y-x)$. For arbitrary $ t, \epsilon > 0$, observe that;

\begin{align*}
\Phi(t) \left|  \Bigl\{ y \in \R^n: \frac{|u(y)L_{x}f(y)|}{u(x)(\|f\|_{ wL_{\Phi}^{u}(\mathbb{R}^n)}+\epsilon)} > t \Bigr\} \right| 
& = \Phi(t) \left|  \Bigl\{ y \in \R^n: \frac{|u(y)f(y-x)|}{u(x)(\|f\|_{wL_{\Phi}^{u}(\mathbb{R}^n)} + \epsilon)} > t \Bigr\} \right| \\
& = \Phi(t) \left|  \Bigl\{ v \in \R^n: \frac{|u(v+x)f(v)|}{u(x)(\|f\|_{wL_{\Phi}^{u}(\mathbb{R}^n)} + \epsilon)} > t \Bigr\} \right| \\
& \leq \Phi(t) \left|  \Bigl\{ v \in\R^n: \frac{|u(v)u(x)f(v)|}{u(x)(\|f\|_{ wL_{\Phi}^{u}(\mathbb{R}^n)} + \epsilon)} > t \Bigr\} \right| \\
& = \Phi(t) \left|  \Bigl\{ v \in \R^n: \frac{|u(v)f(v)|}{(\|f\|_{wL_{\Phi}^{u}(\mathbb{R}^n)} + \epsilon)} > t \Bigr\} \right| \\
& \leq 1
\end{align*}

for $ v:= y-x$. Since $t > 0$ is arbitrary, we have $$\mathop {\sup }\limits_{t > 0} \Phi(t) \left| \left\{ y \in \R^n: \frac{|u(y)L_{x}f(y)|}{u(x)(\|f\|_{wL_{\Phi}^{u}(\mathbb{R}^n)} + \epsilon)} > t \right\} \right|  \leq 1.$$

This shows that $\|L_{x}f\|_{wL_{\Phi}^{u}(\mathbb{R}^n)}  \leq u(x)(\|f\|_{wL_{\Phi}^{u}(\mathbb{R}^n)} + \epsilon)$ for every $ \epsilon > 0$. Hence, we conclude that  $$\|L_{x}f\|_{wL_{\Phi}^{u}(\mathbb{R}^n)} \leq u(x) \|f\|_{wL_{\Phi}^{u}(\mathbb{R}^n)}.$$

{\rm (2)} Let $ f \in wL_{\Phi}^{u}(\mathbb{R}^n)$ and $f\neq0$, then there exist a constant $ C> 0$ (depends on $f$) such that $$\|f\|_{wL_{\Phi}^{u}(\mathbb{R}^n)}\leq C.$$ By Lemma \ref{lemma:2.6} (1), then we have $\|L_{x}f\|_{wL_{\Phi}^{u}(\mathbb{R}^n)} \leq C u(x).$

Now, for every $ t, \epsilon > 0$, we have
\begin{align*}
\Bigl| \Bigl \{ v \in \R^n: \frac{|u(x)f(v)|}{\mathop {\sup }\limits_{ v \in \R^n }u(-v)(\|L_{x}f\|_{wL_{\Phi}^{u}(\mathbb{R}^n)} + \epsilon)} > t \Bigr\}\Bigr| 
& \leq \Bigl| \Bigl\{ v \in \R^n: \frac{|u(x)f(v)|}{u(-v)(\|L_{x}f\|_{wL_{\Phi}^{u}(\mathbb{R}^n)} + \epsilon)} > t \Bigr\} \Bigr|\\
& \leq  \Bigl| \Bigl\{ v \in \R^n: \frac{|u(v+x)f(v)|}{\|L_{x}f\|_{wL_{\Phi}^{u}(\mathbb{R}^n)} + \epsilon} > t \Bigr\}\Bigr|\\
& \leq  \Bigl| \Bigl\{ y \in \R^n: \frac{|u(y)f(y-x)|}{\|L_{x}f\|_{wL_{\Phi}^{u}(\mathbb{R}^n)} + \epsilon} > t \Bigr\} \Bigr|\\
& = \Bigl| \Bigl\{ y \in \R^n: \frac{|u(y)L_{x}f(y)|}{\|L_{x}f\|_{wL_{\Phi}^{u}(\mathbb{R}^n)} + \epsilon} > t \Bigr\} \Bigr|
\end{align*}

for $ y = v+x$. So we obtain, 
$$\Phi(t) \Bigl| \Bigr\{ v \in \R^n: \frac{|u(x)f(v)|}{\mathop {\sup }\limits_{ v \in \R^n}u(-v)(\|L_{x}f\|_{wL_{\Phi}^{u}(\mathbb{R}^n)} + \epsilon)} > t \Bigr\}\Bigr| \leq \Phi(t) \Bigl| \Bigr\{ y \in \R^n: \frac{|u(y)L_{x}f(y)|}{\|L_{x}f\|_{wL_{\Phi}^{u}(\mathbb{R}^n)} + \epsilon} > t \Bigr\} \Bigr|\leq 1.$$

Since $ t > 0$ is arbitrary, we also obtain
  $$\frac{u(x)\|f\|_{wL_{\Phi}(\mathbb{R}^n)}}{\mathop {\sup }\limits_{ v \in \R^n}u(-v)} \leq \|L_{x}f\|_{wL_{\Phi}^{u}(\mathbb{R}^n)},$$ for $\|f\|_{wL_{\Phi}(\mathbb{R}^n)} := \inf \left\{  {b>0:
  	\mathop {\sup }\limits_{t > 0} \Phi(t) \Bigl| \{ x \in \mathbb{R}^n : \frac{|f(x)|}{b} > t \}
  	\Bigr| \leq1}\right\}$.

Choose $C := \mathop {\max } \left\{C_1, \frac{\mathop {\sup }\limits_{ v \in \R^n}u(-v)}{\|f\|_{wL_{\Phi}^{u}(\mathbb{R}^n)}} \right\}$.  Hence, we conclude that $\frac{u(x)}{C} \leq \|L_{x}f\|_{wL_{\Phi}^{u}(\mathbb{R}^n)} \leq C u(x),$ as desired. \qed
\par}

\medskip

Now, we present the sufficient and necessary conditions for the inclusion properties of weighted weak Orlicz spaces $wL_{\Phi_1}^{u_1}(\mathbb{R}^n)$ and $wL_{\Phi_2}^{u_2}(\mathbb{R}^n)$ with respect to Young functions $\Phi_1, \Phi_2$ and weights $u_1, u_2$.

\medskip

\begin{theorem}\label{theorem:2.4}
Let $\Phi_{1}, \Phi_{2}$ be continuous Young functions satisfying the $\bigtriangleup_{2}$ condition such that $\Phi_1 \prec \Phi_2$ and $ u_1, u_2 \in U$. Then the following statements are equivalent:
	
	{\parindent=0cm
		
		{\rm (1)}  $u_1 \preceq u_2$.
		
		{\rm (2)} $wL_{\Phi_2}^{u_2}(\mathbb{R}^n) \subseteq wL_{\Phi_1}^{u_1}(\mathbb{R}^n)$.
		
		{\rm (3)} There exist a constant $C>0$ such that $\|f\|_{wL_{\Phi_1}^{u_1}(\mathbb{R}^n)} \leq C \|f\|_{wL_{\Phi_2}^{u_2}(\mathbb{R}^n)}$, for every $ f \in wL_{\Phi_2}^{u_2}(\mathbb{R}^n)$.
		
		\par}
\end{theorem}

Proof.

Assume that (1) holds. Let  $ f$ be an element of $ wL_{\Phi_2}^{u_2}(\mathbb{R}^n)$. Since $\Phi_1 \prec \Phi_2$ and $u_1 \preceq u_2$, there exists constant $C_1, C_2>0$ such that $\Phi_1(t) \leq \Phi_2(C_1t)$ for all $t > 0$ and $ u_1(x) \leq C_2 u_2(x)$ for every $ x \in \R^n$. Using a similar argument in the proof of Theorem \ref{theorem:2.3} we have $$\|  f \|_{wL_{\Phi_1}^{u_1}(\R^n)} \leq C_1\|  f \|_{wL_{\Phi_2}^{u_1}(\R^n)} \leq C_1C_2 \|  f \|_{wL_{\Phi_2}^{u_2}(\R^n)}.$$ 

As before, we have (2) and (3) are equivalent. It thus remains to show that (3) implies (1).
Assume that (3) holds. By Lemma \ref{lemma:2.6}, we have

$$\frac{u_{1}(x)}{C}\leq \|L_{x}f\|_{wL_{\Phi_1}^{u_1}(\mathbb{R}^n)} \leq C \|L_{x}f\|_{wL_{\Phi_2}^{u_2}(\mathbb{R}^n)}\leq u_2(x),$$ for every $x \in \mathbb{R}^n$. So, we obtain $u_1 \preceq u_2$.\qed

\bigskip

\remark It follows from Theorems \ref{theorem:2.3} and \ref{theorem:2.4} that there cannot be an inclusion relation between $wL_{p_1}^{u_1}(\R^n)$ and $wL_{p_2}^{u_2}(\R^n)$ for distinct values of $p_1$ and $p_2$. In spite of that, for finite measure X we can obtained inclusion relation between $wL_{p_1}^{u_1}(X)$ and $wL_{p_2}^{u_2}(X)$ in the next section.

\section{An additional case}

In the following, we will give sufficient condition for H\"{o}lder's inequality on weighted weak Orlicz spaces which will be used to obtain inclusion relation between $wL_{p_1}^{u_1}(X)$ and $wL_{p_2}^{u_2}(X)$.   

\medskip

\begin{theorem}\label{theorem:3.1} (H\"{o}lder's inequality)
Let $\Phi_1, \Phi_2$, $\Phi_3$ be Young functions and $u_1, u_2, u_3$ be weights such that $\Phi^{-1}_1(t)\Phi^{-1}_2(t) \leq \Phi^{-1}_3(t)$ for every  $t > 0$
and  $u_3(x) \leq u_1(x)u_2(x)$ for every $ x \in X$. If $f_1 \in wL_{\Phi_1}^{u_1}(X)$ and $f_2 \in wL_{\Phi_2}^{u_2}(X)$, then $f_1 f_2 \in wL_{\Phi_3}^{u_3}(X)$ with $$\| f_1 f_2 \|_{wL_{\Phi_3}^{u_3}(X)} \leq 2 \| f_1 \|_{wL_{\Phi_1}^{u_1}(X)} \| f_2 \|_{wL_{\Phi_2}^{u_2}(X)}.$$
\end{theorem}

\noindent{\it Proof}. 

Let $ f_i$ be elements of $ wL_{\Phi_i}^{u_i}(X)$, i =1,2. By Lemma \ref{lemma:2.3}, for every $k \in \mathbb{N}$ we have

\begin{center}
	$\Phi_1(t)\Bigl| \Bigl\{ x \in X : \frac{|u_1(x)f_1(x)|}{(1+\frac{1}{k})\| f_1 \|_{wL_{\Phi_1}^{u_1}(X)} } > t \Bigr\} \Bigr| \leq1$
and
$ \Phi_2(t)\Bigl| \Bigl\{ x \in X : \frac{|u_2(x)f_2(x)|}{(1+\frac{1}{k})\| f_2 \|_{wL_{\Phi_2}^{u_2}(X)} } > t \Bigr\} \Bigr| \leq1$
\end{center}
for every $ t > 0$.

For each $x \in X$ and $k \in \mathbb{N}$, let $$M(x,k) := \max \Bigl(  \Phi_{1}\Bigl(\frac{|u_1(x)f_1(x)|}{(1+\frac{1}{k})\|f_1\|_{wL_{\Phi_1}^{u_1}(X)} }\Bigr) , \Phi_{2}\Bigl(\frac{|u_2(x)f_2(x)|}{(1+\frac{1}{k})\|f_2\|_{wL_{\Phi_2}^{u_2}(X)}}\Bigr)\Bigr).$$

From $\Phi_{i}\Bigl(\frac{|u_i(x)f_i(x)|}{(1+\frac{1}{k})\parallel f_i \parallel_{wL_{\Phi_i}^{u_i}(X)}}\Bigr)\leq M(x,k)$ and Lemma \ref{lemma:1.1} (3), we have $$\frac{|u_i(x)f_i(x)|}{(1+\frac{1}{k})\|f_i\|_{wL_{\Phi_i}^{u_i}(X)}} \leq \Phi_{i}^{-1}\Bigl(\Phi_{i}\Bigl(\frac{|u_i(x)f_i(x)|}{(1+\frac{1}{k})\|f_i\|_{wL_{\Phi_i}^{u_i}(X)}}\Bigr)\Bigr)\leq \Phi_{i}^{-1}(M(x,k)),$$

for $ i = 1,2$.

Hence
$\prod\limits_{i=1}^{2}\frac{|u_i(x)f_i(x)|}{(1+\frac{1}{k})\parallel f_i \parallel_{wL_{\Phi_i}^{u_i}(X)}}\leq \Phi_{1}^{-1}(M(x,k))\Phi_{2}^{-1}(M(x,k)) \leq \Phi_{3}^{-1}(M(x,k))$

and
\begin{align*}
\Phi_{3}\Bigl(\prod\limits_{i=1}^{2}\frac{|u_i(x)f_i(x)|}{(1+\frac{1}{k})\|f_i\|_{wL_{\Phi_i}^{u_i}(X)}}\Bigr) & \leq \Phi_{3}(\Phi_{3}^{-1}(M(x,k)))\\
& \leq M(x,k).
\end{align*}

On the other hand, we have
$M(x,k) \leq\sum\limits_{i=1}^{2} \Phi_{i}\bigl(\frac{|u_i(x)f_i(x)|}{(1+\frac{1}{k})\parallel f_i \parallel_{wL_{\Phi_i}^{u_i}(X)} } \bigr).$

Therefore\\
\begin{align*}
\Phi_{3} (t)\Bigl| \Bigl\{ x \in X : \prod\limits_{i=1}^{2}\frac{\sqrt{|u_3(x)|}|f_i(x)|}{(1+\frac{1}{k})\| f_i \|_{wL_{\Phi_i}^{u_i}(X)}} > t \Bigr\} \Bigr|
&= \Phi_{3} \Bigl(\prod\limits_{i=1}^{2}\frac{\sqrt{|u_3(x)|t_0}|f_i(x)|}{(1+\frac{1}{k})\| f_i \|_{wL_{\Phi_i}^{u_i}(X)}} \Bigr) | \{ x \in X : 1 > t_0 \} |  \\
& \leq \Phi_{3} \Bigl(\prod\limits_{i=1}^{2}\frac{\sqrt{t_0}|u_i(x)f_i(x)|}{(1+\frac{1}{k})\| f_i \|_{wL_{\Phi_i}^{u_i}(X)}}\Bigr) | \{ x \in X : 1 > t_0 \} | \\
& \leq \sum\limits_{i=1}^{2}  \Phi_{i}\Bigl(\frac{\sqrt{t_0}|u_i(x)f_i(x)|}{(1+\frac{1}{k})\| f_i \|_{wL_{\Phi_i}^{u_i}(X)}}\Bigr)| \{ x \in X : 1 > t_0 \}|
\end{align*}

where $ t_0 = \frac{t(1+\frac{1}{k})^{2}\| f_1 \|_{wL_{\Phi_1}^{u_1}(X)} \| f_2 \|_{wL_{\Phi_2}^{u_2}(X)}}{|u_3(x)f_1(x)f_2(x)|}$.\\

Next, we also have\\
\begin{align*}
\Phi_{i}\Bigl(\frac{\sqrt{t_0}|u_i(x)f_i(x)|}{(1+\frac{1}{k})\|f_i\|_{wL_{\Phi_i}^{u_i}(X)}}\Bigr)| \{ x \in X : 1 > t_0 \}| & = \Phi_{i}(t_i) \Bigl| \Bigr\{ x \in X : \Bigl(\frac{|u_i(x)f_i(x)|}{(1+\frac{1}{k})\|f_i\|_{wL_{\Phi_i}^{u_i}(X)}} \Bigr)^2 > t_i^2 \Bigr\} \Bigr| \\
& = \Phi_{i}(t_i) \Bigl| \Bigl\{ x \in X : \frac{|u_i(x)f_i(x)|}{(1+\frac{1}{k})\|f_i\|_{wL_{\Phi_i}^{u_i}(X)}}  > t_i \Bigr\} \Bigr| \\
& \leq 1
\end{align*}

where $ t_i = \frac{\sqrt{t_0}|u_i(x)f_i(x)|}{(1+\frac{1}{k})\parallel f_i \parallel_{wL_{\Phi_i}^{u_i}(X)}}$, for $i =1, 2$.\\

So, we obtain $\Phi_{3} (t)\Bigl| \Bigl\{ x \in X : \prod\limits_{i=1}^{2}\frac{\sqrt{|u_3(x)|}|f_i(x)|}{(1+\frac{1}{k})\| f_i \|_{wL_{\Phi_i}^{u_i}(X)}} > t \Bigr\} \Bigr| \leq 2$.

On the other hand, we have
\begin{align*}
\Phi_{3} (t)\Bigl| \Bigl\{ x \in X : \prod\limits_{i=1}^{2}\frac{\sqrt{|u_3(x)|}|f_i(x)|}{\sqrt{2}(1+\frac{1}{k})\| f_i \|_{wL_{\Phi_i}^{u_i}(X)}} > t \Bigr\} \Bigr|
 \leq & \mathop {\sup }\limits_{t > 0}\Phi_{3} (t)\Bigl| \Bigl\{ x \in X : \prod\limits_{i=1}^{2}\frac{\sqrt{|u_3(x)|}|f_i(x)|}{(1+\frac{1}{k})\| f_i \|_{wL_{\Phi_i}^{u_i}(X)}} > 2t \Bigr\} \Bigr| \\
 = & \mathop {\sup }\limits_{s > 0}\Phi_{3} \Bigl(\frac{s}{2}\Bigr)\Bigl| \Bigl\{ x \in X : \prod\limits_{i=1}^{2}\frac{\sqrt{|u_3(x)|}|f_i(x)|}{(1+\frac{1}{k})\| f_i \|_{wL_{\Phi_i}^{u_i}(X)}} > s \Bigr\} \Bigr| \\
  \leq & \mathop {\sup }\limits_{s > 0}\frac{1}{2}\Phi_{3} (s)\Bigl| \Bigl\{ x \in X : \prod\limits_{i=1}^{2}\frac{\sqrt{|u_3(x)|}|f_i(x)|}{(1+\frac{1}{k})\| f_i \|_{wL_{\Phi_i}^{u_i}(X)}} > s \Bigr\} \Bigr|\\
 \leq & 1.
\end{align*}

Since $ t > 0$ is an arbitrary positive real number, we get

$$\mathop {\sup }\limits_{t > 0}\Phi_{3} (t)\Bigl| \Bigl\{ x \in X : \frac{|u_3(x)f_1(x)f_2(x)|}{2(1+\frac{1}{k})^{2}\| f_1 \|_{wL_{\Phi_1}^{u_1}(X)}\| f_2 \|_{wL_{\Phi_2}^{u_2}(X)}} > t \Bigr\} \Bigr| \leq 1.$$

This shows that $\|f_1f_2\|_{wL_{\Phi_3}^{u_3}(X)} \leq 2(1+\frac{1}{k})^{2}\|f_1\|_{wL_{\Phi_1}^{u_1}(X)}\|f_2\|_{wL_{\Phi_2}^{u_2}(X)}$ and this is true for every $k \in \mathbb{N}$. We can conclude that   $\|f_1f_2\|_{wL_{\Phi_3}^{u_3}(X)} \leq
2 \| f_1 \|_{wL_{\Phi_1}^{u_1}(X)}\|  f_2 \|_{wL_{\Phi_2}^{u_2}(X)}$, as desired. \qed

\bigskip

\begin{corollary}\label{corollary:3.2}
	Let $X:=B(a,r_0) \subseteq \mathbb{R}^n$ for some $a\in\mathbb{R}^n$ and $r_{0} > 0$.
	If $\Phi_1, \Phi_2$ are
	two Young functions, $u_1, u_2$ are weights and there are a Young function $\Phi$ and a weight $0 < u(x) \leq 1$ for every $ x \in X$ such that
	$
	\Phi^{-1}_1(t)\Phi^{-1}(t) \leq \Phi^{-1}_2(t)
	$
	for every $ t \geq 0$ and $ u_1(x) \leq u(x)u_2(x)$ for every $ x\in X$, then $$wL_{\Phi_1}^{u_1}(X)\subseteq wL_{\Phi_2}^{u_2}(X)$$ with
	$
	\| f \|_{wL_{\Phi_2}^{u_2}(X)} \leq \frac{2}{\Phi^{-1}(\frac{1}
		{|B(a,r_{0})|})} \| f \|_{wL_{\Phi_1} ^{u_1}(X)}
	$
	for $ f \in wL_{\Phi_1}^{u_1}(X)$.
\end{corollary}

\noindent{\it Proof}. Since $ 0 < u(x) \leq 1$ for every $ x \in X$, we have $ \|f\|_{wL_{\Phi_2}^{u_2}(X)} \leq \|\frac{f}{u}\|_{wL_{\Phi_2}^{u_2}(X)}$. 
Let $f \in wL_{\Phi_1}^{u_1}(X)$, by Theorem \ref{theorem:3.1} and choosing $g:=\chi_{B(a,r_0)}$,
we obtain
\begin{align*}
\| f\|_{wL_{\Phi_2}^{u_2}(X)} = & \| f \chi_{B(a,r_0)} \|_{wL_{\Phi_2}^{u_2}(X)} \\
\leq & \left\|\frac{f\chi_{B(a,r_0)}}{u} \right\|_{wL_{\Phi_2}^{u_2}(X)} \\  
\leq & 2 \left\|\frac{\chi_{B(a,r_0)}}{u}
\right\|_{wL_{\Phi}^u(X)} \| f \|_{wL_{\Phi_1}^{u_1}(X)} \\
= &\frac{2}{\Phi^{-1}(\frac{1}{|B(a,r_{0})|})} \| f \|_{wL_{\Phi_1}^{u_1}(X)}.
\end{align*}

This shows that $wL_{\Phi_1}^{u_1}(X) \subseteq wL_{\Phi_2}^{u_2}(X)$.

\medskip

We shall now discuss the inclusion properties of weighted weak Lebesgue spaces $wL_{p_1}^{u_1}(X)$ and $wL_{p_2}^{u_2}(X)$ with respect to distinct values of $p_1$ and $p_2$ as well as $u_1$ and $u_2$.

\medskip

\begin{corollary}
	Let $X:=B(a,r_0)$ for some $a\in\mathbb{R}^n$ and $r_0>0$. If $ 1 \leq p_{2} <
	p_{1} < \infty$ and $u_1(x) \leq u_2(x)$ for every $ x \in X$, then $wL_{p_{1}}^{u_1}(X) \subseteq wL_{p_{2}}^{u_2}(X)$.
\end{corollary}

\noindent{\it Proof}.
Let $\Phi_1(t):= t^{p_1}, \Phi_2(t):= t^{p_2}$, and $\Phi(t):=
t^{\frac{p_1 p_2}{p_1- p_2}}$ ($t\ge 0$). Since $ 1 \leq p_{2} < p_{1} < \infty$,
we have $\frac{p_1 p_2}{p_1- p_2} > 1$. Thus, $\Phi_1,\
\Phi_2$, and $\Phi$ are Young functions. Now, define $u(x) = \frac{u_1(x)}{u_2(x)}$ for every $x \in X$. Observe that, using the definition
of $\Phi^{-1}$ and Lemma \ref{lemma:1.1}, we have
$$
\Phi^{-1}_{1}(t) =t^{\frac{1}{p_1}},\ \Phi^{-1}_{2}(t) =t^{\frac{1}{p_2}},\
{\rm and}\ \Phi^{-1}(t) =t^{\frac{p_1- p_2}{p_1p_2}}.
$$
Moreover, $\Phi^{-1}_{1}(t)\Phi^{-1}(t) = t^{\frac{1}{p_1}}t^{\frac{p_1-
		p_2}{p_1p_2}} = t^{\frac{1}{p_2}} = \Phi^{-1}_{2}(t)$ and $ u_1(x) = \frac{u_1(x)}{u_2(x)}u_2(x) = u(x)u_2(x)$. So it follows
from Corollary \ref{corollary:3.2} that $\| f \|_{wL_{p_2}^{u_2}(X)} \leq \frac{2}{\Phi^{-1}(
	\frac{1}{|B(a,r_0)|})}  \| f \|_{wL_{p_1}^{u_1}(X)}$, and therefore we can conclude that $wL_{p_{1}}^{u_1}(X)
\subseteq wL_{p_{2}}^{u_2}(X)$. \qed

\section{Concluding Remarks}

We have shown the sufficient and necessary conditions for the inclusion relation between weighted weak Lebesgue spaces and between weighted weak Orlicz
spaces. The inclusion properties of weighted weak Orlicz spaces generalize the inclusion properties of weak Orlicz spaces in \cite{Masta1}. In the proof of the inclusion property we used the norm of characteristic function on $\R^n$ and estimate the norm of the translation functions in $\R^n$. As a corollary of Theorem \ref{theorem:2.4} and Theorem 2.22 in \cite{Alen}, we have that the inclusion properties of weighted weak Orlicz spaces are equivalent with the inclusion properties of weighted Orlicz spaces.

\end{document}